# Decomposing Intraday Dependence in Currency Markets: Evidence from the AUD/USD Spot Market

Jonathan A. Batten [ab*], Craig A. Ellis [c], Warren P. Hogan [d]

[a] Graduate School of Management, Macquarie University, CBD Office, Level 3, 51-57 Pitt St, Sydney Australia 2000

[b] College of Business Administration, Seoul National University, Sillim-Dong, Kwanak-ku, 151-742, South Korea.

[c] School of Economics and Finance, University of Western Sydney, Locked Bag 1797, Penrith South DC, NSW 1797, Australia

[d] School of Finance and Economics, University of Technology Sydney, Haymarket 2007 NSW, Australia

[*] Corresponding author. Jonathan A. Batten, Graduate School of Management, Macquarie University, CBD Office, Level 3, 51-57 Pitt St., Sydney Australia 2000. Phone: ++61-2-8274-8344, Fax: ++61-2-8274-8370, Email: jonathan.batten@mgsm.edu.au

# Decomposing Intraday Dependence in Currency Markets:
# Evidence from the AUD/USD Spot Market

**Abstract**

The local Hurst exponent, a measure employed to detect the presence of dependence in a time series, may also be used to investigate the source of intraday variation observed in the returns in foreign exchange markets. Given that changes in the local Hurst exponent may be due to either a time-varying range, or standard deviation, or both of these simultaneously, values for the range, standard deviation and local Hurst exponent are recorded and analyzed separately. To illustrate this approach, a high-frequency data set of the spot Australian dollar/U.S. dollar provides evidence of the returns distribution across the 24-hour trading 'day', with time-varying dependence and volatility clearly aligning with the opening and closing of markets. This variation is attributed to the effects of liquidity and the price-discovery actions of dealers.

**Decomposing Intraday Dependence in Currency Markets:**

**Evidence from the AUD/USD Spot Market**

1. **Introduction**

The availability of high-frequency or tick data of spot foreign exchange (FX) prices over the last decade has expanded the range of empirical investigation possible in a foreign exchange markets. For example in the financial economics literature, Muller *et al.* [1], Goodhart and Demos [2], Goodhart and Figliuoli [3] and later Bollerslev and Domowitz [4] and Bollerslev and Melvin [5] focus on quote arrivals (frequency) and the size of the bid-ask spread, which they find varies across the trading day, with higher spreads and volatility at the beginning and end of trading. From a time-series modelling perspective, differences in liquidity and price availability of markets organised around groups of dealers - who possess differing degrees of private information - ensures that prices cannot immediately incorporate all private information in individual trades [6]. Price discovery by these traders may therefore lead to time-series that display statistical properties consistent with dependent processes.

Studies investigating the statistical properties of financial series (Cheung [7]; Batten and Ellis [8]; Van De Gucht, Dekimpe, and Kwok [9]; Opong, Mulholland, Fox, and Farahmand [10]; Barkoulas, Labys, and Onochie [11]) identify the presence of non-linear dependence, which is a departure from the fair game, or martingale property of asset returns under Fama's [12] Efficient Market Hypothesis. In the econophysics literature, recent studies focusing on the long-range dependent properties of stock indices by Grau-Carles [13], Costa and Vasconcelos [14], Matos *et al* [15] and Cajueiro and Tabak [16,17,18] also describe varying levels of long-range dependence.

The implications of dependent processes, evident from low and high order autocorrelation structures in the data are of particular concern for the volatility based pricing models (such as option pricing models) typically used in financial markets. Low order correlations, which tend to exhibit hyperbolic decay, may be associated with short-range memory effects, while long-range memory effects have been linked to the presence of fractal structures.

Despite some studies investigating these issues in the major traded currencies such as the Euro, Japanese yen, or English pound quoted against the U.S. dollar, there is little information available on the microstructure and statistical properties of trading on the spot Australian dollar against the U.S. dollar (AUD/USD). Extending the work of Batten and Ellis [19] who found weak evidence of positive dependence in the daily spot AUD/USD, this study provides evidence on the intraday behaviour of volatility and links this to the time-varying nature of dependence evident in the series. Measured using the statistical techniques of Hurst [20] and Mandelbrot and Wallis [21], the unique feature of this study is that we decompose the measure of dependence into its underlying components to provide an insight into the cause of the observed variation in volatility and dependence over the 24-hour trading day. The approach differs from recent studies that seek to identify dependence by various methods (see Refs. [13],[15]) and from those which adopt a rolling sample approach (see Refs. [14],[16],[17],[22]), the latter of which generally fail to account for the asymptotic behaviour of the Hurst statistic [23,24]. First we track the patterns and distributions of price quotes, spreads and returns across the trading day and week. Second we describe the nature and the form of price dependence in the markets. The approach adopted is to investigate the statistical relationship between quote returns across the trading day and week using the rescaled range technique, which is then decomposed to identify the source of the observed time variation. We believe that this is the first study to investigate dependence in this manner.

The data employed is time stamped spot AUD/USD price quotes from banks contributing to the Reuters "FX=" page from Friday 5[th] May 2000 to Thursday 15[th] June 2000[1]. The AUD/USD is the sixth most actively traded currency (after the U.S. dollar, euro, yen, Swiss franc and the Canadian dollar) and is traded 24-hours a day with most trading occurring outside the Australian time zone, in the United Kingdom and the United States. Although seemingly a short length of physical time, the data set of approximately 30,000 observations is in fact significantly larger than many recently published studies investigating dependence. Furthermore, the large number of observations examined is more than statistically sufficient to enable us to first measure dependence using the Hurst-Mandelbrot-Wallis [20,21] method and then decompose it into its component parts: the range and standard deviation of intraday returns. While there are some studies using larger high-frequency data sets than employed herein (see for example Ref. [2]), none investigate the specific issues undertaken in our study. In this way the analysis undertaken is pioneering for highlighting further opportunities to be explored and the results offer hypotheses, which, however tentative, may be subject to further work. The proprietary nature of the information from which the current data set is obtained has largely prevented subsequent and other detailed analyses in other markets.

The remainder of the paper is structured as follows: In the next section, the method of calculating dependence and returns is described. In Section 3 the Australian dollar FX market is briefly described, then the characteristics of the quotes in the spot AUD/USD market are established. Then the results from an analysis of variance on the return properties across the trading day and week are presented. The final section allows for some concluding remarks.

---

[1] This sample period avoided the public holiday in Japan on the 5 May (Children's Day) though it included the Monday 29 May holiday in the US (Memorial Day) and UK (Bank Holiday), and the holiday in Australia on the 12 June (Queen's Birthday).

## 2. Research method

*2.1. Defining and measuring long-range dependence.*

Financial economics has attempted to develop fundamental linear models for understanding the relationship between the various variables that are observable in complex economic systems. These linear models have the advantage of providing unique and single solutions. Improved statistical techniques have enabled researchers to accommodate dynamic and non-linear qualities to describe variations on these simple linear models. These variations propose some series are non-linear dynamic possessing certain qualities including:

i. the past memory of events influence the present, such that innovations possess both a transitory component, the effect of which slowly dissipates over time, plus a permanent component; and

ii. chaotic systems where complex non-linear relations may be indistinguishable from random behaviour when using standard statistical tests.

Non-linear dynamic behaviour may follow either a deterministic or stochastic path. A usual feature of deterministic non-linear behaviour is long-range dependence between time-series increments. Time-series that follow a standard Brownian motion (random walk) typically exhibit autocorrelation functions that tend quickly to zero. By contrast, the autocorrelation functions of long-range dependent series tend to decay only slowly to zero. Rather than being independent, the increments of long-range dependent time-series exhibit non-negligible statistical dependence even when they are far away from each other. More formally a series can be defined as a long-range dependent process if for a finite constant $C$, the autocorrelation function, $\rho(k)$ decays as

$$\rho(k) \sim Ck^{2H-2} \qquad (1)$$

as $k \to \infty$. The variable $H$ (Hurst exponent) is a real value $0 \leq H \leq 1$ that describes the sign and magnitude of dependence in the series. For series exhibiting positive long-range dependence or persistence, the value of the exponent is $H>0.5$. Anti-persistent series are conversely characterised by exponent values $H<0.5$. Independent series are finally characterised by $H = 0.5$. Of importance to note is that the statistical definition of long-range dependence in (1) has little to do with physical time (ie. the length of time over which observations are recorded), but instead relies upon the number of historic observations over which a statistically significant correlation can be determined. Thus, distinguishing between the length of physical time over which observations are recorded and the number of observations themselves allows deterministic non-linear behaviour to be examined even over relatively short time frames.

Mathematically one can easily demonstrate that a large number of observations are more important than the length of physical time over which they are recorded. This may best be conceived in terms of the asymptotic nature of the Hurst statistic noted in [23,24]. The critical point is that the length of time over which observations are recorded should only be considered important in terms of the interpretation of findings. To provide economic meaning, we conveniently bundle the estimated Hurst exponent (and its components) into hourly averages (see Table 3, and later illustrated in Figure 3) to demonstrate its time varying nature. This variation conveniently aligns with the opening and closing of markets, which enables us to provide an economic interpretation of our findings in terms of price discovery. The price discovery actions of market participants over the trading day have been well proven in microstructure studies [2,3] employing large high-frequency data sets over longer time periods.

Long-range dependence may be measured using statistical techniques based on range analysis [20,21,25], spectral regression [26], variance estimators [27], maximum likelihood techniques [28,29], and wavelets [30]. The statistical method employed for measuring long-range dependence in this study is based on the classical rescaled adjusted range [20,21]. Being non-parametric in nature the classical rescaled adjusted range has the advantages of being robust to series with significant third and fourth moments [31] and to stochastic processes with infinite variance [32], both of which are common features of financial time-series [33].

The classical rescaled adjusted range $(R/\sigma)_n$ is calculated as

$$(R/\sigma)_n = (1/\sigma_n) \left[ \underset{1 \leq k \leq n}{\text{Max}} \sum_{j=1}^{k}(X_j - \overline{X_n}) - \underset{1 \leq k \leq n}{\text{Min}} \sum_{j=1}^{k}(X_j - \overline{X_n}) \right] \qquad (2)$$

where $\overline{X}_n$ is the sample mean $(1/n)\Sigma_j X_j$ and $\sigma_n$ is the series standard deviation

$$\sigma_n = \left[ 1/n \sum_{j=1}^{n}(X_j - \overline{X_n})^2 \right]^{0.5} \qquad (3)$$

For a given series of length $N$, the procedure first requires that equations (2) and (3) be estimated over several subseries of length $N \geq n$. Using Ordinary Least Squares regression, a global Hurst exponent, $H$ is then estimated. The exponent value is global in the sense that the value of $H$ describes the mean level of dependence for the entire series. In order to capture the time-varying nature of intraday dependence this study employs a local measure of the Hurst exponent, $h_n$

$$h_n = \frac{\log(R/\sigma)_n}{\log n} \quad (4)$$

Estimates of the local Hurst exponent are calculated for ($N$–$n$+1) times overlapping subseries of length $n$, with $n$ having a set value. In this study $n$ is arbitrarily set to either 10-quotes or 20-quotes. The procedure in effect creates a time-series of exponent values, the change in whose value can be measured over time. Given that changes in the local Hurst exponent may be due to a time-varying range, $R_n$ standard deviation, $\sigma_n$, or both of these simultaneously, values for the range, standard deviation and local Hurst exponent are separately recorded for the series.

One criticism of the classical rescaled adjusted range methodology by Lo [25] has been the general failure by researchers to provide for comparison, empirical results generated from a random series. Following a method adopted by Ancel, Ambrose and Griffiths [34], the AUD/USD returns series under investigation is scrambled 1,000 times and values of $R_n$, $\sigma_n$, and $h_n$ recorded at each iteration. Bootstrapped values of the local Hurst exponent, $\hat{h}_n$ for $n$ equals 10 and 20-quotes are subsequently calculated and the moments of the distribution of the exponent estimates used to construct upper and lower confidence intervals. The mean and standard deviation of the local Hurst exponent when $n$ = 10 is 0.4753 and 0.0163 respectively, and for $n$ = 20 is 0.5058 and 0.0197 respectively. Comparison of the observed means, $\bar{h}_n$ across a 24-hour trading day to their bootstrap expected values, $E(\bar{h}_n)$ is undertaken using a one-sample $Z$-test to compute a confidence interval and to perform the following hypothesis test of the mean:

$$H_0: \quad \bar{h}_n = E(\bar{h}_n) \qquad \forall\, n = 10, 20$$

$$H_1: \quad \bar{h}_n \neq E(\bar{h}_n) \qquad \forall\, n = 10, 20 \qquad (5)$$

*2.2    Estimating Returns on Tick Data Quotes*

The approach adopted in this study in dealing with the irregularity of time intervals between quotes is to consider transaction time separately to clock time (see Ref. [35]). In line with the argument that the time between trades (which is a measure of trading activity) affects market price behaviour, we maintain the clock time of the series and measure the actual change in time between quotes, $\tau$ based on a fraction of a 24-hour day. By way of example 12:00 would be 0.5000 (or 12/24) and 13:00 would be 0.5417 (or 13/24), such that the clock time between 13:00 and 12:00, $\tau = 0.0417$.

There are also a number of different approaches to estimating the return on tick data quotes (see Refs. [36],[37]). We report hourly returns adjusted for the $\tau$ interval, measured as a fraction of the 24-hour trading day expressed in seconds, between 2 quotes multiplied by 360 seconds in the hour (to improve the reporting of the statistics). Therefore given the $\tau$ interval between consecutive price quotes expressed in seconds, then the Adjusted Return (AR) is

$$AR = [\log(Q_i) - \log(Q_{(i+\tau)})] * \frac{360}{\tau} \qquad (6)$$

The method is consistent with Andersen and Bollerslev [37] who argue that the logarithmic transformation induces normality in the series by minimising kurtosis.

*2.3    Time variation in returns and local Hurst statistics*

To provide economic meaning to the stream of Hurst statistics generated across the sample period, we conveniently bundle the estimated Hurst exponent (and its components) into hourly averages. Analysis of variance (ANOVA) is then used to investigate and model the relationship between a response variable and one or more independent variables. The method extends the two-sample *t*-test for the equality of two population means to a more general null hypothesis of comparing the equality of more than two means, versus them not all being equal. That is the null hypothesis of similarity between groups, k where the groups are hourly intervals across the trading day is

$$H_0: \mu_1 = \mu_2 = \ldots \mu_k \qquad \forall\, k$$

$$H_1: \mu_1 \neq \mu_2 \neq \ldots \mu_k \qquad \forall\, k \qquad (7)$$

The *F*-test *p*-value is employed to indicate the degree of significant of differences among the means.

**3.    Data and sample**

The best estimate of the scale and scope of the Australian dollar FX market is provided by the Bank for International Settlements [38] triennial survey of foreign exchange and derivatives market activity. A key feature of these surveys is the extent to which trade in Australian dollars is now conducted in financial centres outside Australia. Turnover in Australian dollars against the U.S. dollar was US$47 billion per day, with US$27 billion occurring in Australian time zone trading. This reflects an increase of 14 per cent over April 1998, an outcome boosted by the relocation by a number of global players of their Asian time zone foreign exchange business to Australia. The market is highly concentrated with the ten largest dealers accounting for 76 per cent of total market turnover, slightly down on the level of 80 per cent in 1998. The extent that price formation now occurs in offshore centres, usually outside the time zone for normal Australian trading, is also evident from the high frequency data presented in this study.

The spot AUD/USD quotes used in this study were provided by Reuters Australia Ltd. and comprise nearly 30,000 time stamped (Day:Hour:Minute:Second) spot Australian dollar price quotes from contributing banks based on Greenwich mean time (GMT) from Friday $5^{th}$ May 9:49:11AM GMT 2000 (19:49 Sydney time) to Thursday $15^{th}$ June 0:56:6AM GMT, 2000 (10:56 Sydney time). This comprises 42 calendar days. The price quotes were used to update the FX= page on the Reuters Terminal (RT) system. Foreign exchange spot trading generally begins the week by 23:00 GMT on Sunday (9:00am Sydney time) and ends the week usually by 22.00pm GMT Friday night (5.00pm New York). The latest Friday night quote was on the 5/5/2000 at 21:50). Therefore key times across the trading day are the opening of markets in Sydney (23:00 GMT), Hong Kong and Singapore (24:00 GMT), Tokyo (1:00 GMT), Europe (8:00 GMT), London (9:00 GMT), and New York (14:00 GMT). The London and New York markets then close at 17:00 GMT and 22:00 GMT respectively.

Details of the number of quotes generated on a weekday basis are provided in Table 1. Excluding the Sunday morning (early morning trading in Sydney, Australia) with 348 quotes over the sample, the average number of quotes of 1032 is about one quarter the quantity generated in a major currency (such as the 4,500 reported by [35] on the USD/DMK). Most quotes were received on Tuesday and the least quotes were received on the Friday. Batten and Hogan [39] provide further economic discussion of the data set and the intraday actions of informed and uninformed traders.

*** *(Insert Table 1 about here)* ***

The time intervals of the quotes are not regular and the time periods were later converted to fractions of a 24-hour day based on GMT. The time stamp also enabled each observation to be coded with a proxy for each day of the week (from 1 to 7 with 1 = Monday and Sunday = 7). The quote frequency also varied over the trading day. This variation is clearly visible in Figure 1, which plots the average number of quotes received per hour of the sample period. The four peaks in this diagram at 1:00, 6:00 and 12:00 and 14:00 GMT represent the opening of the Tokyo market (1:00 GMT), late afternoon Asia and early morning European trading (6:00 GMT), and then close of Asian trading and the period before lunchtime in London (12:00 GMT). The peak frequency (around 14:00 GMT) was in late afternoon trading in London and mid-morning trading in U.S. and Canadian markets.

*** *(Insert Figure 1 about here)* ***

Table 1 also reports the mean arrival time and the coefficient of variation between quotes, $\mu_\tau$ of the arrival time for each hour in the day. The $F$-statistic based on a one-way ANOVA suggests significant statistical variation with the period of greatest quote arrivals being the period from 12:00 to 15:00 GMT which is the period of overlap between London afternoon trading and early morning U.S. trading. It is interesting to note that the last few hours of trading have few and infrequent quotes with a high coefficient of variation of arrival. This improves, as would be expected, with the opening of markets in Asia, Europe and the United States.

*** (Insert Table 2 about here) ***

4. Results

The ANOVA of the mean adjusted returns, AR and the local Hurst coefficients based on a 10 and 20-quote interval, across a 24-hour trading day, are reported in Table 2. The Hurst coefficient provides a measure of the degree of dependence in the market and may also proxy as a measure of efficiency or departure from a random walk, which has a Hurst coefficient of $H = 0.5$. The $F$-statistic of the ANOVA of the mean hourly returns, suggest that the trading day is associated with statistically significant differences (at the 99.9% level of confidence) in the mean hourly returns and variance. This result was also supported by a non-parametric Kruskal-Wallis test of differences in the median returns across the trading day. The Kruska-Wallis $H$-statistic of 3146 with 23 degrees of freedom was significant at the 99.9% level of confidence. The volatility appears to be greatest at 21:00 GMT with spikes in volatility at 3:00, 10:00 and 14:00 GMT. GMT 21:00 is 4:00PM in the east coast of the U.S. and Canada and is associated with the close of North American trading, a decrease in quote arrival, 3:00 GMT is 12:00 noon in Tokyo and the therefore coincides with period prior to Tokyo lunchtime, 10:00 GMT is the period following the opening of the London market and 14:00 GMT is the opening of the east coast U.S. market and the period following the London lunch-time. This last period is the period of minimum quote arrival time, and maximum quote arrival time.

The results suggest three distinct sources of volatility in the AUD/USD spot market. The first is liquidity-induced volatility associated with the absence of quotes as occurs in after-hours trading (such as the period between the close of New York and the opening of the Australian markets). The second is volatility associated with price-discovery at the start of trading. This occurs at the opening of markets and is evident from trading at the opening of markets in Australia, Asia and London markets. In these instances volatility arises from the information seeking activities of uninformed traders. These traders would be expected to congregate when trading costs are least [35]. The third form of volatility is associated with price-discovery at the close of markets as dealers balance or close their order-books. The spike in volatility at 21:00 GMT provides a clear example. These relationships are evident from Figure 2, which plots the mean and the variance of returns. These graphs display the typical 'U-shaped' volatility pattern that has been observed in other FX markets.

*** *(Insert Figure 2 and 3 about here)* ***

Table 2 provides the mean, $\bar{h}_n$ and standard deviation, $\sigma(\hat{h}_n)$ of the local Hurst coefficient calculated over a 10 and 20-quote interval. The *F*-statistic of both the 10 and 20-quote interval local Hurst suggests statistically significant differences in the mean coefficients across the trading day. This is consistent with time-varying dependence as well as the previously noted time variation in volatility. This is clearly evidenced in Figure 3, which plots the 10 and 20-quote Hurst coefficient across the trading day. The Figure also includes the estimated Hurst coefficient based upon a 1000 iteration randomizing of the return series. In line with Lo (1991) a confidence interval was also estimated based on a 1000 iteration randomizing of the series. The *Z*-test of acceptance of the null hypothesis that the $\bar{h}_n = E(\bar{h}_n)$ for $n = 10$ and $n = 20$ is rejected in favour of the actual values not equally the expected values of 0.4753 and 0.5058 respectively. The Z-statistic for $n = 10$ at the 95% confidence interval is 102.41, significant at the 99.99% level, and for $n = 20$ is 164.63, also significant at the 99.99% level.

A further implication of this result is that evidence of significant dependence across the trading day, suggested by observed values of $h_n$ being greater than the $\bar{h}_n$ in the return series, cannot be accounted for by transience in the classical rescaled adjusted range alone.

*** *(Insert Table 3 about here)* ***

The source of variation in the local Hurst statistic across the trading day is provided in Table 3, which provides a deconstruction of the exponent based on the mean of the range and the standard deviation for 10 and 20-quote intervals across the trading day. The first four columns report the raw range, $R_n$ and standard deviations, $\sigma_n$ with their respective logarithms based on a 10-quote interval. The next four columns report the same information based on a 20-quote interval. The *F*-statistic of an ANOVA of differences in the means across the trading day clearly provides evidence of statistically significant differences in the hourly means based on a 10 and 20-quote interval.

The change in the range, $DR_n$ and standard deviation, $D\sigma_n$ logarithm values between consecutive hours are then reported in the next four columns for the 10-quote and 20-quote interval. It will be apparent from these percentages that the source of the variation in the local Hurst exponent is due to changes in variance for both the 10 and 20-quote construction of the exponent. The next four columns report the percentage change in the range and the standard deviation for the local Hurst exponent constructed from 10 and 20-quotes, $DR_{10-20}$ and $D\sigma_{10-20}$ respectively. The source of variation between the local Hurst exponent constructed from 10 and 20-quotes is the clearly the range; as the quote interval increases the Hurst exponent will approach its asymptotic value. The last two columns report the local Hurst exponent calculated as the difference in the logarithms of the range and standard deviation, divided by the logarithm of the quote interval, being either 10 or 20-quotes. Overall the results suggest that the time-varying nature of dependence in the series investigated is due to the time-varying nature of volatility, which is due to the effects of changes in liquidity and the price discovery actions of dealers.

## 5.   Conclusions

The study investigates the source of intraday variation commonly observed in the returns in foreign exchange markets when high-frequency data sets are employed. To illustrate our approach, the differences in returns, volatility and dependence in quotes on the spot AUD/USD exchange rate are compared across the 24-hour trading day. The sample covered tick data from the 5[th] May 2000 to the 15 June 2000 and comprised approximately 30,000 observations. Over this period, quote arrival is clearly –and not surprisingly- linked to the opening and closing of key markets in Asia, Europe and the U.S., with the most dense quote period being late afternoon London as the U.S. markets open. The frequency of quote arrival also varied with the longest time between quotes during the late afternoon trading in the U.S. just prior to the opening of the Australian markets.

Attention is drawn to 'U-shaped' distribution of volatility over the trading day and the consequent time-varying nature of dependence induced by this variation in volatility. The statistically significant difference in the mean local Hurst exponent across the 24-hour trading day is not attributable to transience in the estimation of the exponent. Decomposing the local exponent into its 'range' and 'standard deviation' components reveals that changes in the mean standard deviation across the trading day are largely responsible for the variation in the local exponent. In discussion of this finding, three underlying sources of variation in AUD/USD spot volatility are suggested: volatility associated with the absence of quotes, and volatility associated with price-discovery at the start of trading and at the close of markets.

Table 1
Cross-Tabulation of Spot AUD/USD Quotes by Hour and Weekday
(Sample Period 5 May 2000 to 15 June 2000)

| GMT | Monday | Tuesday | Wednesday | Thursday | Friday | Sunday | Total | Mean Arrival Time in Minutes $\mu_\tau$ | Coefficient of Variation of Arrival Time |
|---|---|---|---|---|---|---|---|---|---|
| 0 | 187 | 229 | 190 | 160 | 75 | 0 | 841 | 1.84 | 1.36 |
| 1 | 246 | 177 | 228 | 190 | 72 | 40 | 953 | 1.68 | 1.44 |
| 2 | 178 | 121 | 161 | 140 | 95 | 12 | 707 | 2.20 | 1.61 |
| 3 | 157 | 64 | 79 | 52 | 38 | 0 | 390 | 4.06 | 1.42 |
| 4 | 158 | 129 | 136 | 85 | 61 | 0 | 569 | 3.05 | 1.54 |
| 5 | 176 | 198 | 212 | 174 | 98 | 0 | 858 | 1.94 | 1.45 |
| 6 | 344 | 441 | 484 | 423 | 288 | 0 | 1980 | 0.83 | 1.39 |
| 7 | 356 | 441 | 442 | 333 | 264 | 0 | 1836 | 0.88 | 1.23 |
| 8 | 335 | 340 | 357 | 331 | 263 | 0 | 1626 | 0.99 | 1.27 |
| 9 | 273 | 365 | 292 | 232 | 244 | 0 | 1406 | 1.16 | 1.24 |
| 10 | 253 | 373 | 339 | 232 | 222 | 0 | 1419 | 1.18 | 1.44 |
| 11 | 257 | 419 | 398 | 345 | 327 | 0 | 1746 | 0.98 | 1.22 |
| 12 | 393 | 566 | 572 | 474 | 481 | 0 | 2486 | 0.68 | 1.11 |
| 13 | 389 | 547 | 549 | 458 | 437 | 0 | 2380 | 0.70 | 1.10 |
| 14 | 420 | 483 | 522 | 593 | 1117 | 0 | 3135 | 0.85 | 1.75 |
| 15 | 354 | 445 | 469 | 304 | 461 | 0 | 2033 | 0.82 | 1.14 |
| 16 | 323 | 359 | 286 | 249 | 290 | 0 | 1507 | 1.10 | 1.09 |
| 17 | 176 | 212 | 225 | 143 | 199 | 0 | 955 | 1.74 | 1.01 |
| 18 | 87 | 127 | 115 | 72 | 109 | 0 | 510 | 3.23 | 0.80 |
| 19 | 74 | 74 | 61 | 44 | 37 | 20 | 310 | 5.43 | 1.26 |
| 20 | 40 | 42 | 55 | 33 | 22 | 25 | 217 | 8.47 | 1.05 |
| 21 | 72 | 84 | 73 | 34 | 3 | 64 | 330 | 6.28 | 1.56 |
| 22 | 118 | 118 | 108 | 73 | 0 | 74 | 491 | 3.34 | 1.29 |
| 23 | 195 | 314 | 141 | 128 | 0 | 113 | 891 | 1.96 | 1.87 |
| Total | 5561 | 6668 | 6494 | 5302 | 5203 | 348 | 29576 | | |

The Table reports the average number of quotes (N) received across the 24-hour trading day based on Greenwich Mean Time (GMT). The column labeled $\mu_\tau$ is the mean arrival time in seconds between consecutive quotes ($Q_i$ and $Q_{i+\tau}$). The *F*-statistic of equality of the means in arrival times is 289.92, a *p*-value = 0.000 based on a one-way analysis of variance (ANOVA). The *p*-value indicates that there is significant evidence at the 99.9% level of differences among the means of arrival time.

Table 2
ANOVA of Spot AUD/USD Returns and the Local Hurst Exponent Across the Trading Day
(Sample Period 5 May 2000 to 15 June 2000)

| GMT | $\mu_R$ | $\sigma_R$ | $\bar{h}_{10}$ | $\sigma(\hat{h}_{10})$ | $\bar{h}_{20}$ | $\sigma(\hat{h}_{20})$ | Bootstrap $\bar{h}_{10}$ | Bootstrap $\sigma(\hat{h}_{10})$ | Bootstrap $\bar{h}_{20}$ | Bootstrap $\sigma(\hat{h}_{20})$ |
|---|---|---|---|---|---|---|---|---|---|---|
| 0 | 4.592 | 6.242 | 0.4866 | 0.0614 | 0.5035 | 0.0446 | 0.4749 | 0.0183 | 0.5057 | 0.0189 |
| 1 | 4.194 | 6.042 | 0.4892 | 0.0755 | 0.5187 | 0.0613 | 0.4747 | 0.0157 | 0.5077 | 0.0196 |
| 2 | 5.510 | 8.874 | 0.4805 | 0.0758 | 0.5361 | 0.0609 | 0.4776 | 0.0158 | 0.5059 | 0.0193 |
| 3 | 10.153 | 14.454 | 0.4855 | 0.0816 | 0.5403 | 0.0713 | 0.4748 | 0.0166 | 0.5054 | 0.0198 |
| 4 | 7.626 | 11.732 | 0.4960 | 0.0737 | 0.5326 | 0.0703 | 0.4751 | 0.0161 | 0.5057 | 0.0192 |
| 5 | 4.851 | 7.052 | 0.4899 | 0.0714 | 0.5171 | 0.0693 | 0.4748 | 0.0159 | 0.5063 | 0.0212 |
| 6 | 2.063 | 2.872 | 0.4896 | 0.0764 | 0.5289 | 0.0717 | 0.4741 | 0.0170 | 0.5050 | 0.0199 |
| 7 | 2.196 | 2.699 | 0.4832 | 0.0805 | 0.5227 | 0.0718 | 0.4748 | 0.0165 | 0.5057 | 0.0203 |
| 8 | 2.464 | 3.130 | 0.4862 | 0.0776 | 0.5263 | 0.0687 | 0.4788 | 0.0169 | 0.5060 | 0.0186 |
| 9 | 2.899 | 3.593 | 0.4888 | 0.0708 | 0.5308 | 0.0675 | 0.4743 | 0.0156 | 0.5055 | 0.0197 |
| 10 | 2.948 | 4.244 | 0.4880 | 0.0803 | 0.5363 | 0.0715 | 0.4746 | 0.0163 | 0.5047 | 0.0198 |
| 11 | 2.447 | 2.993 | 0.4730 | 0.0801 | 0.5064 | 0.0690 | 0.4750 | 0.0173 | 0.5048 | 0.0197 |
| 12 | 1.692 | 1.882 | 0.4804 | 0.0777 | 0.5183 | 0.0718 | 0.4756 | 0.0143 | 0.5060 | 0.0191 |
| 13 | 1.756 | 1.940 | 0.4847 | 0.0794 | 0.5200 | 0.0699 | 0.4736 | 0.0167 | 0.5060 | 0.0192 |
| 14 | 2.119 | 3.700 | 0.4840 | 0.0803 | 0.5186 | 0.0667 | 0.4747 | 0.0167 | 0.5063 | 0.0200 |
| 15 | 2.046 | 2.339 | 0.4776 | 0.0817 | 0.5137 | 0.0719 | 0.4756 | 0.0160 | 0.5072 | 0.0196 |
| 16 | 2.758 | 3.006 | 0.4832 | 0.0848 | 0.5221 | 0.0766 | 0.4757 | 0.0153 | 0.5047 | 0.0195 |
| 17 | 4.344 | 4.406 | 0.5009 | 0.0897 | 0.5622 | 0.0758 | 0.4737 | 0.0156 | 0.5054 | 0.0193 |
| 18 | 8.074 | 6.442 | 0.4675 | 0.1190 | 0.5382 | 0.0959 | 0.4759 | 0.0167 | 0.5055 | 0.0197 |
| 19 | 13.582 | 17.051 | 0.4956 | 0.0968 | 0.5275 | 0.0779 | 0.4729 | 0.0164 | 0.5072 | 0.0201 |
| 20 | 21.162 | 22.226 | 0.4910 | 0.0903 | 0.5673 | 0.0762 | 0.4783 | 0.0164 | 0.5056 | 0.0199 |
| 21 | 15.713 | 24.591 | 0.4855 | 0.0755 | 0.5414 | 0.0745 | 0.4750 | 0.0161 | 0.5059 | 0.0196 |
| 22 | 8.361 | 10.796 | 0.4951 | 0.0813 | 0.5453 | 0.0633 | 0.4758 | 0.0162 | 0.5060 | 0.0201 |
| 23 | 4.893 | 9.160 | 0.4956 | 0.0797 | 0.5421 | 0.0676 | 0.4773 | 0.0180 | 0.5055 | 0.0214 |
| | | | | | | | | | | |
| F-Statistic | 289.93 | | 8.31 | | 39.84 | | | | | |
| p-value | 0.000 | | 0.000 | | 0.000 | | | | | |
| $E(\bar{h}_n)$ | | | | | | | 0.4753 | | 0.5058 | |
| $\sigma(\hat{h}_n)$ | | | | | | | | 0.0164 | | 0.0197 |

The Table reports the average hourly returns and the Hurst coefficient across the 24-hour trading day, and a set of bootstrap local Hurst coefficients based on a 1000 iteration randomizing of the time-series following Lo (1991). The column labeled $\mu_R$ is the mean hourly return, and the column labeled $\sigma_R$ is the mean hourly standard deviation between consecutive quotes ($Q_i$ and $Q_{i+\tau}$). The mean hourly local Hurst coefficient calculated on a 10-quote range ($\bar{h}_{10}$), and a 20-quote range ($\bar{h}_{20}$), while the adjoining columns report the standard deviation of the 10 and 20-quote means ($\sigma(\hat{h}_{10})$ and $\sigma(\hat{h}_{20})$ respectively). The F-statistics report a one-way analysis of variance (ANOVA) that tests for equality of means across the 24-hours. The F-statistic p-value indicates that there is significant evidence at the 99.9% level of differences among the means of returns and the Hurst coefficient. The last four columns report the mean and standard deviation of the simulated (bootstrap) series for each hour and for the overall sample.

Table 3
Source of Variation in the Local Hurst Exponent (Sample Period 5 May 2000 to 15 June 2000)

| GMT | Range $R_{10}$ | log $R_{10}$ | Std Dev $\sigma_{10}$ | log $\sigma_{10}$ | Range $R_{20}$ | log $R_{20}$ | Std Dev $\sigma_{20}$ | log $\sigma_{20}$ | Change $DR_{10}$ | Change $D\sigma_{10}$ | Change $DR_{20}$ | Change $D\sigma_{20}$ | $DR_{10-20}$ | $D\sigma_{10-20}$ | Local Hurst $h_{10}$ | Local Hurst $h_{20}$ |
|---|---|---|---|---|---|---|---|---|---|---|---|---|---|---|---|---|
| 0 | 963.3 | 2.984 | 321.02 | 2.507 | 2042.9 | 3.310 | 468.96 | 2.671 | | | | | 10.94 | 6.57 | 0.4772 | 0.4912 |
| 1 | 114.3 | 2.058 | 37.89 | 1.579 | 183.3 | 2.263 | 41.65 | 1.620 | -31.03 | -88.20 | -31.63 | -39.37 | 9.97 | 2.60 | 0.4795 | 0.4946 |
| 2 | 47.3 | 1.675 | 15.62 | 1.194 | 82.2 | 1.915 | 18.03 | 1.256 | -18.62 | -58.78 | -15.39 | -22.45 | 14.33 | 5.22 | 0.4812 | 0.5064 |
| 3 | 29.1 | 1.464 | 9.20 | 0.964 | 43.4 | 1.638 | 8.57 | 0.933 | -12.60 | -41.10 | -14.49 | -25.72 | 11.86 | -3.20 | 0.5001 | 0.5415 |
| 4 | 29.2 | 1.465 | 9.20 | 0.964 | 52.7 | 1.722 | 10.66 | 1.028 | 0.10 | 0.00 | 5.15 | 10.16 | 17.50 | 6.64 | 0.5016 | 0.5335 |
| 5 | 17.4 | 1.241 | 5.54 | 0.744 | 32.7 | 1.515 | 6.94 | 0.841 | -15.34 | -39.78 | -12.04 | -18.14 | 22.09 | 13.16 | 0.4970 | 0.5174 |
| 6 | 6.3 | 0.799 | 2.01 | 0.303 | 12.4 | 1.093 | 2.49 | 0.396 | -35.57 | -63.72 | -27.81 | -52.91 | 36.79 | 30.67 | 0.4961 | 0.5359 |
| 7 | 5.9 | 0.771 | 1.88 | 0.274 | 10.1 | 1.004 | 2.02 | 0.305 | -3.56 | -6.47 | -8.15 | -22.93 | 30.29 | 11.38 | 0.4967 | 0.5372 |
| 8 | 7.2 | 0.857 | 2.31 | 0.364 | 12.9 | 1.111 | 2.59 | 0.413 | 11.22 | 22.87 | 10.58 | 35.35 | 29.54 | 13.67 | 0.4937 | 0.5360 |
| 9 | 8.3 | 0.919 | 2.64 | 0.422 | 14.4 | 1.158 | 2.88 | 0.459 | 7.20 | 14.29 | 4.30 | 11.15 | 26.04 | 8.96 | 0.4975 | 0.5372 |
| 10 | 8.9 | 0.949 | 2.86 | 0.456 | 16.3 | 1.212 | 3.15 | 0.498 | 3.30 | 8.33 | 4.65 | 8.47 | 27.68 | 9.19 | 0.4930 | 0.5487 |
| 11 | 7.1 | 0.851 | 2.33 | 0.367 | 12.7 | 1.104 | 2.67 | 0.427 | -10.34 | -18.53 | -8.94 | -14.41 | 29.67 | 16.10 | 0.4839 | 0.5206 |
| 12 | 4.5 | 0.653 | 1.48 | 0.170 | 8.0 | 0.903 | 1.65 | 0.218 | -23.27 | -36.48 | -18.18 | -49.01 | 38.25 | 27.74 | 0.4830 | 0.5270 |
| 13 | 4.6 | 0.663 | 1.48 | 0.170 | 7.8 | 0.892 | 1.60 | 0.204 | 1.46 | 0.00 | -1.22 | -6.14 | 34.60 | 19.89 | 0.4925 | 0.5288 |
| 14 | 12.0 | 1.079 | 3.97 | 0.599 | 23.8 | 1.377 | 5.26 | 0.721 | 62.83 | 168.24 | 54.31 | 253.22 | 27.56 | 20.41 | 0.4804 | 0.5039 |
| 15 | 5.5 | 0.740 | 1.77 | 0.248 | 9.1 | 0.959 | 1.89 | 0.277 | -31.40 | -55.42 | -30.33 | -61.66 | 29.54 | 11.49 | 0.4924 | 0.5246 |
| 16 | 7.0 | 0.845 | 2.21 | 0.344 | 11.7 | 1.068 | 2.35 | 0.371 | 14.15 | 24.86 | 11.38 | 34.22 | 26.40 | 7.75 | 0.5007 | 0.5358 |
| 17 | 10.0 | 1.000 | 3.00 | 0.477 | 18.0 | 1.255 | 3.17 | 0.501 | 18.33 | 35.75 | 17.51 | 35.03 | 25.53 | 5.02 | 0.5229 | 0.5797 |
| 18 | 13.3 | 1.124 | 4.35 | 0.639 | 23.8 | 1.377 | 4.52 | 0.655 | 12.39 | 45.00 | 9.66 | 30.75 | 22.49 | 2.61 | 0.4854 | 0.5545 |
| 19 | 153.0 | 2.185 | 50.40 | 1.702 | 214.5 | 2.331 | 60.34 | 1.781 | 94.39 | 1058.62 | 69.36 | 171.79 | 6.72 | 4.59 | 0.4823 | 0.4234 |
| 20 | 167.3 | 2.224 | 54.64 | 1.738 | 307.8 | 2.488 | 76.53 | 1.884 | 1.78 | 8.41 | 6.73 | 5.80 | 11.91 | 8.42 | 0.4860 | 0.4646 |
| 21 | 163.2 | 2.213 | 53.52 | 1.729 | 368.0 | 2.566 | 88.69 | 1.948 | -0.48 | -2.05 | 3.12 | 3.40 | 15.96 | 12.69 | 0.4842 | 0.4750 |
| 22 | 32.0 | 1.505 | 9.94 | 0.997 | 148.2 | 2.171 | 31.26 | 1.495 | -31.98 | -81.43 | -15.39 | -23.25 | 44.23 | 49.89 | 0.5078 | 0.5195 |
| 23 | 16.9 | 1.228 | 5.33 | 0.727 | 38.5 | 1.586 | 7.60 | 0.881 | -18.42 | -46.38 | -26.97 | -41.08 | 29.12 | 21.20 | 0.5012 | 0.5416 |
| F-statistic | 408.91 | | 408.42 | | 1222.69 | | 1151.53 | | | | | | | | | |
| p-value | 0.000 | | 0.000 | | 0.000 | | 0.000 | | | | | | | | | |

The Table reports details of the construction of the local Hurst based on the mean of the Range and the Standard Deviation for 10 and 20-quote intervals across the trading day.

Figure 1
Plot of the Total Number of Quotes Across the Trading Day
(Sample Period 5 May 2000 to 15 June 2000)

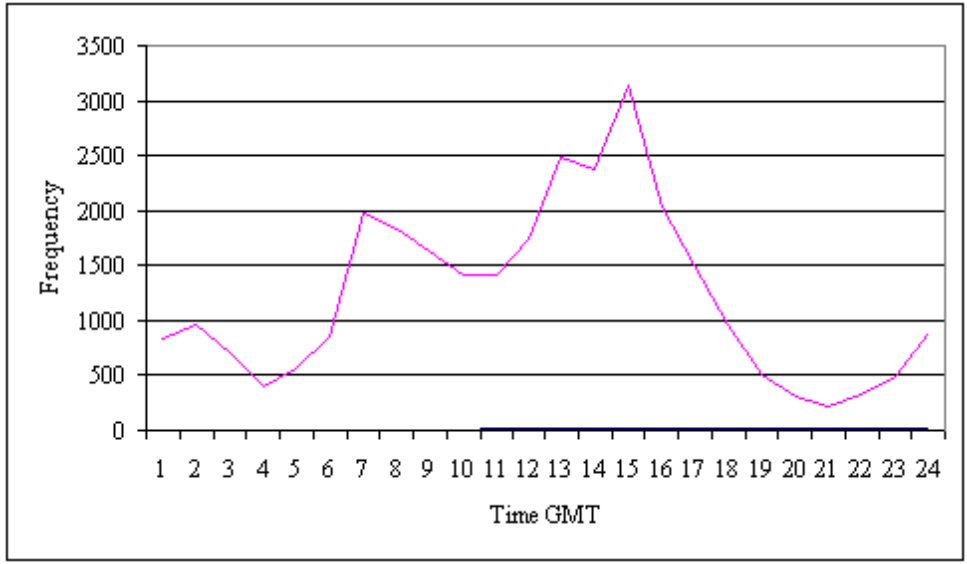

The Figure plots the average frequency of quote arrivals across the 24 hour trading day.

Figure 2
Plot of the Adjusted AUD/USD Return Mean and Volatility Across the Trading Day
(5 May 2000 to 15 June 2000)

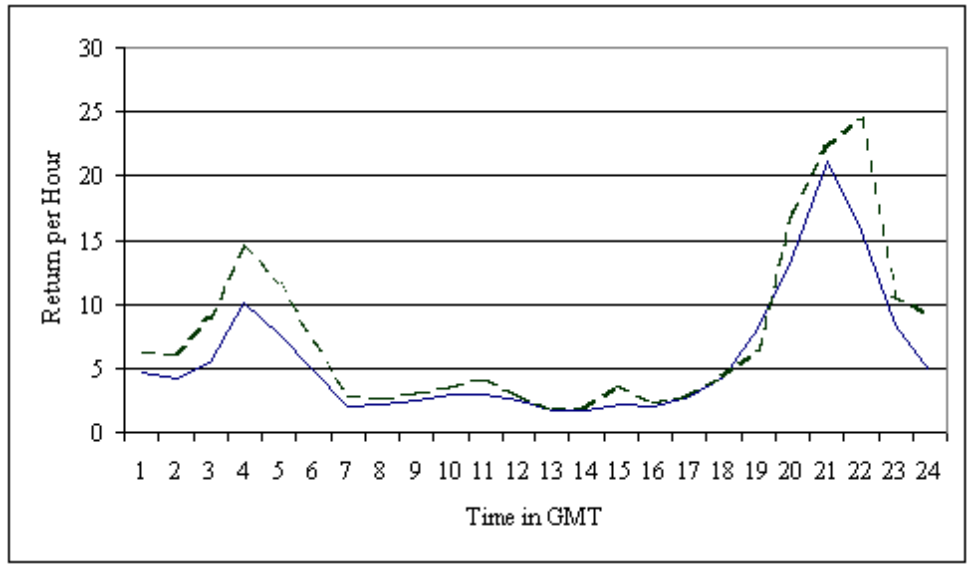

The Figure plots the average hourly returns and the volatility of returns across the 24-hour trading day. The volatility is the dashed line and the mean is the connected line. The return is the adjusted hourly log difference of the midrate between two consecutive price quotes.

Figure 3
Plot of the Mean of the Local Hurst Exponent
Based on a 10 and 20 Quote Interval Across the Trading Day
(5 May 2000 to 15 June 2000)

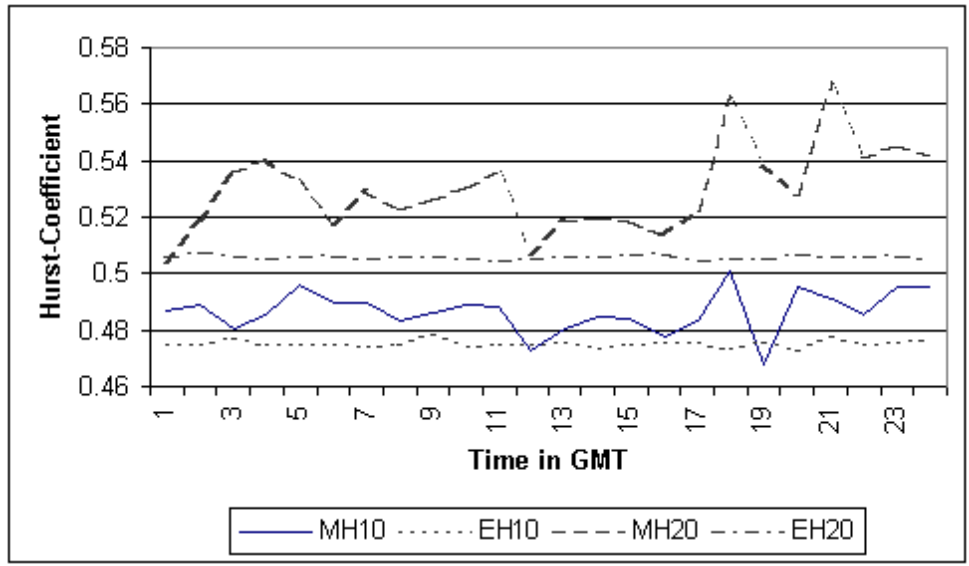

The Figure plots the mean of Hurst (MH) coefficient across the 24-hour trading day. The mean hourly Hurst coefficient is calculated on a 10-quote range (MH10), and a 20-quote range (MH20). Also plotted is the simulated Hurst for the 10 and 20-quote range based on a randomizing of the time series. The top dashed line (MH20) is greater than the straight line (MH10) and both series exceed the expected Hurst coefficient for 46 of the 48 hour values.